\documentclass[12pt,reqno]{amsart}

\usepackage{amssymb}

\usepackage{eucal}

\input{diagrams}

\font\sss=cmss8

\def\cA{{\mathcal A}}
\def\cB{{\mathcal B}}
\def\cC{{\mathcal C}}

\def\cH{{\mathcal H}}
\def\cI{{\mathcal I}}

\def\cM{{\mathcal M}}
\def\cN{{\mathcal N}}
\def\cO{{\mathcal O}}

\def\sA{\mbox{\sf A}}

\def\sD{\mbox{\sf D}}

\def\ast{{\textstyle *}}

\def\D{\sD}

\def\Dsmall{\mbox{\sss D}}

\def\End{\operatorname{End}}
\def\Ext{\operatorname{Ext}}

\def\Hom{\operatorname{Hom}}
\def\id{\operatorname{id}}

\def\LTensor{\stackrel{\operatorname{L}}{\otimes}}
\def\Mod{\mbox{\sf Mod}}
\def\Modsmall{\mbox{\sss Mod}}

\def\QCoh{\mbox{\sf QCoh}}
\def\QCohsmall{\mbox{\sss QCoh}}

\def\Rf{\operatorname{R}\!f_{\ast}}
\def\RGamma{\operatorname{R}\!\Gamma}

\def\RSHom{\operatorname{R}\!\cH\mbox{\it om}}

\def\Spec{\operatorname{Spec}}

\numberwithin{equation}{part}

\newtheorem{Lemma}{Lemma}[section]
\newtheorem{Theorem}[Lemma]{Theorem}
\newtheorem{Proposition}[Lemma]{Proposition}
\newtheorem{Corollary}[Lemma]{Corollary}

\theoremstyle{definition}

\def\AR{Aus\-lan\-der-Rei\-ten}

\def\ARs{\AR\ se\-quen\-ce}
\def\ARth{\AR\ the\-o\-ry}

\def\cat{\sA}  

\def\Homcat{\Hom}
\def\Endcat{\End}
\def\Extcat{\Ext}

\def\otimessch{\otimes}

\def\HomQCoh{\Hom}
\def\EndQCoh{\End}
\def\ExtQCoh{\Ext}

\def\sch{X}  

\def\leftobj{\cA}    
\def\midobj{\cB}     
\def\midobjprime{\cB^{\prime}}     
\def\rightobj{\cC}   

\def\injobj{\cI}     
\def\testobj{\cM}    

\def\testcomp{\cN}   

\def\leftobjmor{a}   
\def\leftobjmorprime{\alpha}  
\def\midobjmor{b}    
\def\rightobjmor{c}  

\def\injobjmor{i}    
\def\injobjmorprime{\iota}    
\def\testobjmor{m}   
\def\sigmaobjmor{\sigma}      

\def\DX{\sD(\QCoh\,\sch)}
\def\DsmallX{\mbox{\sss D}(\QCohsmall\,\sch)}

\begin{document}

\title[Auslander-Reiten sequences]
{Auslander-Reiten sequences on schemes}

\author{Peter J\o rgensen}
\address{Department of Pure Mathematics, University of Leeds,
Leeds LS2 9JT, United Kingdom}
\email{popjoerg@hotmail.com, www.geocities.com/popjoerg}


\keywords{Auslander-Reiten theory, duality, smooth projective
scheme, quasi-coherent sheaf, smooth projective curve}

\subjclass[2000]{14F05, 16G70}

\begin{abstract} 
Let $\sch$ be a smooth projective scheme of dimension $d \geq 1$ over
the field $k$, and let $\rightobj$ be an indecomposable coherent sheaf
on $\sch$.  Then there is an \ARs\ in the category of
quasi-coherent sheaves on $\sch$,
\[
  0 \rightarrow (\Sigma^{d-1}\rightobj) \otimessch \omega
    \longrightarrow \midobj
    \longrightarrow \rightobj
    \rightarrow 0.
\]
Here $\Sigma^{d-1}\rightobj$ is the $(d-1)$'st syzygy in a minimal
injective resolution of $\rightobj$, and $\omega$ is the dualizing
sheaf of $\sch$.
\end{abstract}

\maketitle

\setcounter{section}{-1}
\section{Introduction}
\label{sec:introduction}

This note shows that \ARs s frequently exist in categories of
quasi-coherent sheaves on schemes.  More precisely, let $\sch$ be a
smooth projective scheme of dimension $d \geq 1$ over the field $k$,
and let $\rightobj$ be an indecomposable coherent sheaf on $\sch$.
Then by theorem \ref{thm:ARs}, there is an \ARs\ in the category
of quasi-coherent sheaves on $\sch$,
\begin{equation}
\label{equ:ARs}
  0 \rightarrow \leftobj 
    \longrightarrow \midobj
    \longrightarrow \rightobj
    \rightarrow 0.
\end{equation}
Moreover, $\leftobj$ can be computed: It is $(\Sigma^{d-1}\rightobj)
\otimessch \omega$, where $\Sigma^{d-1}\rightobj$ is the
$(d-1)$'st syzygy in a minimal injective resolution of $\rightobj$ in
the category of quasi-coherent sheaves, and $\omega$ is the dualizing
sheaf of $\sch$.

The sheaves $\leftobj$ and $\midobj$ are not in general coherent, but
only quasi-coherent.  This is analogous to ring theory: If $C$ is a
finitely presented non-projective $R$-module with local endomorphism
ring, then by \cite[thm.\ 4]{ARing} there is an \ARs\ in the category
of all $R$-modules,
\[
  0 \rightarrow A \longrightarrow B \longrightarrow C \rightarrow 0,
\] 
but $A$ and $B$ are not in general finitely presented.

However, note that if $\sch$ is a curve, then $d = 1$, and then
$\Sigma^{d-1}\rightobj$ is just $\rightobj$ which {\em is} coherent,
so in this case, $\leftobj$ and $\midobj$ are coherent.  So if $\sch$
is a curve, then I recover the result known from
\cite{ReitenVandenBergh} that the category of coherent sheaves on
$\sch$ has \ARs s; see corollary \ref{cor:ARs}.

The proof that the \ARs\ \eqref{equ:ARs} exists is based on
proposition \ref{pro:ARs} which uses another form of duality than the
classical Auslander-Reiten formula to get \ARs s.  This may be of
independent interest.

\section{\ARs s and duality}
\label{sec:higher_duality}

\begin{Proposition}
\label{pro:ARs}
Let $\cat$ be a $k$-linear abelian category with enough injectives,
over the field $k$.  Let $\rightobj$ be an object with local
endomorphism ring.  Let $\leftobj$ be another object for which there
is a natural equivalence
\begin{equation}
\label{equ:equivalence}
  \Homcat(\rightobj,-)^{\prime} \simeq \Extcat^d(-,\leftobj)
\end{equation}
for some $d \geq 1$, where the prime denotes dualization with respect
to $k$.

Then there is a short exact sequence
\[
  0 \rightarrow \Sigma^{d-1}\leftobj 
    \longrightarrow \midobj
    \stackrel{\midobjmor}{\longrightarrow} \rightobj
    \rightarrow 0
\]
where $\midobjmor$ is right almost split.  Here $\Sigma^{d-1}\leftobj$
is the $(d-1)$'st syzygy in an injective resolution of $\leftobj$.
\end{Proposition}

\begin{proof}
Denote by $J$ the Jacobson radical of the endomorphism ring
$\Homcat(\rightobj,\rightobj)$, and pick a non-zero linear map
$\epsilon$ in $\Homcat(\rightobj,\rightobj)^{\prime}$ which vanishes
on $J$.  By the natural equivalence \eqref{equ:equivalence} the map
$\epsilon$ corresponds to a non-zero element $e$ in
$\Extcat^d(\rightobj,\leftobj)$.

Let $0 \rightarrow \leftobj \longrightarrow \injobj^0
\longrightarrow \cdots \longrightarrow \injobj^{d-1} 
\stackrel{\injobjmor}{\longrightarrow} \Sigma^d\leftobj \rightarrow 0$
be terms number $0$ to $d-1$ of an injective resolution of $\leftobj$,
augmented by $\leftobj$ and the $d$'th syzygy $\Sigma^d\leftobj$.

Pick a morphism $\rightobj \stackrel{\rightobjmor}{\longrightarrow}
\Sigma^d\leftobj$ which represents the element $e$ in
$\Extcat^d(\rightobj,\leftobj)$, and form the commutative diagram
\[
  \begin{diagram}[labelstyle=\scriptstyle]
    0 & \rTo & \Sigma^{d-1}\leftobj & \rTo & \midobj 
      & \rTo^{\midobjmor} & \rightobj & \rTo & 0 \\
      & & \dEq & & \dTo & & \dTo_{\rightobjmor} & & \\
    0 & \rTo & \Sigma^{d-1}\leftobj & \rTo & \injobj^{d-1} 
      & \rTo_{\injobjmor} & \Sigma^d\leftobj & \rTo & 0 \lefteqn{,} \\ 
  \end{diagram}
\]
where the lower row comes from the injective resolution of $\leftobj$,
and where the right hand square is formed by pull back; cf.\
\cite[prop.\ VIII.4.2]{MacLane}.

The morphism $\midobjmor$ is not split, for otherwise it would be easy
to see that $\rightobjmor$ factored through $\injobjmor$, but then
$\rightobjmor$ would represent $0$ in $\Extcat^d(\rightobj,\leftobj)$,
contradicting that $\rightobjmor$ represents $e$ which is non-zero.

Now let $\testobj \stackrel{\testobjmor}{\longrightarrow} \rightobj$
be a morphism which is not a split epimorphism.  Then 
$
\begin{diagram}[labelstyle=\scriptstyle]
\Homcat(\rightobj,\testobj)
& \rTo^{\Homcat(\rightobj,\testobjmor)} &
\Homcat(\rightobj,\rightobj)
\end{diagram}
$
is not surjective, for $\id_{\rightobj}$ is not in the image.  Now,
$\Homcat(\rightobj,\testobjmor)$ is clearly a homomorphism of
right-modules over the ring $\Homcat(\rightobj,\rightobj)$.  The
target of this homomorphism is $\Homcat(\rightobj,\rightobj)$ itself
which is a local ring.  Since $\Homcat(\rightobj,\testobjmor)$ is not
surjective, its image must be contained in $J$, which is (among other
things) the unique maximal proper right-submodule of
$\Homcat(\rightobj,\rightobj)$.  

Hence the linear map $\epsilon$ vanishes on the image of
$\Homcat(\rightobj,\testobjmor)$, so
$\Homcat(\rightobj,\testobjmor)^{\prime}$ sends $\epsilon$ to $0$.  By
the equivalence \eqref{equ:equivalence} this means that
$\Extcat^d(\testobjmor,\leftobj)$ sends $e$ to $0$.  Now, $e$ is
represented by $\rightobj \stackrel{\rightobjmor}{\longrightarrow}
\Sigma^d\leftobj$, so $\Extcat^d(\testobjmor,\leftobj)(e)$ is
represented by $\testobj
\stackrel{\rightobjmor\testobjmor}{\longrightarrow}
\Sigma^d\leftobj$, so $\Extcat^d(\testobjmor,\leftobj)(e) = 0$ implies
that 
$\testobj \stackrel{\rightobjmor \testobjmor}{\longrightarrow}
\Sigma^d \leftobj$ 
factors through 
$\injobj^{d-1} \stackrel{\injobjmor}{\longrightarrow} \Sigma^d \leftobj$.
So there is a commutative diagram
\[
  \begin{diagram}[labelstyle=\scriptstyle]
    \testobj & & & & \\
    & \rdTo(4,2)^{\testobjmor} \rdTo(2,4) & & & \\
    & & \midobj & \rTo^{\midobjmor} & \rightobj \\
    & & \dTo & & \dTo_{\rightobjmor} \\
    & & \injobj^{d-1} & \rTo_{\injobjmor} & \Sigma^d\leftobj \lefteqn{.} \\ 
  \end{diagram}
\]
But the square is a pullback square, so the diagram can be completed
with a morphism $\testobj \longrightarrow \midobj$, so
$\testobjmor$ factors through $\midobjmor$.
\end{proof}

\section{A lemma on endomorphism rings}
\label{sec:endomorphism_rings}

\begin{Lemma}
\label{lem:End}
Let $\cat$ be an abelian category with a short exact sequence 
\[
  0 \rightarrow \leftobj 
    \stackrel{\leftobjmor}{\longrightarrow} \injobj 
    \stackrel{\injobjmor}{\longrightarrow} \Sigma\leftobj
    \rightarrow 0, 
\]
where $\leftobj \stackrel{\leftobjmor}{\longrightarrow} \injobj$ is an
injective envelope and where $\Extcat^1(\injobj,\leftobj) = 0$.

If $\leftobj$ has local endomorphism ring, then so does 
$\Sigma \leftobj$. 
\end{Lemma}

\begin{proof}
Consider the commutative diagram
\[
  \begin{diagram}[labelstyle=\scriptstyle]
    0 & \rTo & \leftobj & \rTo^{\leftobjmor} & \injobj 
      & \rTo^{\injobjmor} & \Sigma\leftobj & \rTo & 0 \\
    & & \dTo^{\leftobjmorprime} 
      & & \dTo_{\injobjmorprime} 
      & & \dTo_{\sigmaobjmor} & & \\
    0 & \rTo & \leftobj & \rTo_{\leftobjmor} & \injobj 
      & \rTo_{\injobjmor} & \Sigma\leftobj & \rTo & 0 \lefteqn{.} \\
  \end{diagram}
\]
It is clear that any morphism $\leftobjmorprime$ gives rise to such a
diagram.  However, so does any morphism $\sigmaobjmor$, for when
$\sigmaobjmor$ is given, consider $\injobj \stackrel{\sigmaobjmor
\injobjmor}{\longrightarrow} \Sigma\leftobj$ which represents an element in
$\Extcat^1(\injobj,\leftobj)$.  Since this $\Ext$ is $0$, the morphism
$\injobj \stackrel{\sigmaobjmor\injobjmor}{\longrightarrow}
\Sigma\leftobj$ must factor through $\injobj 
\stackrel{\injobjmor}{\longrightarrow} \Sigma\leftobj$.  This gives
the morphism $\injobj \stackrel{\injobjmorprime}{\longrightarrow}
\injobj$, and the morphism  $\leftobj
\stackrel{\leftobjmorprime}{\longrightarrow} \leftobj$ and hence the
diagram follows.

Observe that in the diagram, if $\leftobjmorprime$ is an isomorphism,
then so is $\injobjmorprime$ because $\leftobj
\stackrel{\injobjmor}{\longrightarrow} \injobj$ is an injective
envelope, and hence, so is $\sigmaobjmor$.  

To show that $\Endcat(\Sigma\leftobj)$ is local, I must show that if
$\sigmaobjmor$ is a non-invertible element, then $\id_{\Sigma\leftobj}
- \sigmaobjmor$ is invertible.  That $\sigmaobjmor$ is non-invertible
means that it is not an isomorphism.  Embed $\sigmaobjmor$ in a
diagram as above.  Then $\leftobjmorprime$ is not an isomorphism, for
otherwise $\sigmaobjmor$ would be an isomorphism, as observed
above. So $\leftobjmorprime$ is a non-invertible element of
$\Endcat(\leftobj)$.  But then $\id_{\leftobj} - \leftobjmorprime$ is
invertible because $\Endcat(\leftobj)$ is local.  That is,
$\id_{\leftobj} - \leftobjmorprime$ is an isomorphism.  But there is
also a commutative diagram
\[
  \begin{diagram}[labelstyle=\scriptstyle]
    0 & \rTo & \leftobj & \rTo^{\leftobjmor} & \injobj 
      & \rTo^{\injobjmor} & \Sigma\leftobj & \rTo & 0 \\
    & & \dTo^{\id_{\leftobj} - \leftobjmorprime} 
      & & \dTo_{\id_{\injobj} - \injobjmorprime} 
      & & \dTo_{\id_{\Sigma\leftobj} - \sigmaobjmor} & & \\
    0 & \rTo & \leftobj & \rTo_{\leftobjmor} & \injobj 
      & \rTo_{\injobjmor} & \Sigma\leftobj & \rTo & 0 \lefteqn{,} \\
  \end{diagram}
\]
and as observed above, when $\id_{\leftobj} - \leftobjmorprime$ is an
isomorphism, so is $\id_{\Sigma\leftobj} - \sigmaobjmor$.  That is,
$\id_{\Sigma\leftobj} - \sigmaobjmor$ is an invertible element of
$\Endcat(\Sigma\leftobj)$.
\end{proof}

\section{Schemes}
\label{sec:schemes}

The following lemma uses the dualizing sheaf of a projective scheme,
as described for instance in \cite[sec.\ III.7]{Hartshorne}.

\begin{Lemma}
\label{lem:duality}
Let $\sch$ be a projective scheme of dimension $d$ with Gorenstein
singularities, over the field $k$.  Let $\omega$ be the dualizing
sheaf.

Let $\rightobj$ be a coherent sheaf which has a bounded resolution of
locally free coherent sheaves.  Then there are natural equivalences
\[
  \Ext_{\QCohsmall\,\sch}^i(\rightobj,-)^{\prime}
  \simeq \Ext_{\QCohsmall\,\sch}^{d-i}(-,\rightobj \otimessch \omega),
\]
where $\QCoh(\sch)$ is the category of quasi-coherent sheaves on
$\sch$.
\end{Lemma}

\begin{proof}
Since $\sch$ is projective over $k$, there is a projective morphism
$\sch \stackrel{f}{\longrightarrow} \Spec(k)$.  This is certainly a
separated morphism of quasi-compact separated schemes, so according to
\cite[exam.\ 4.2]{NeemanDuality} the derived global section functor of
$\sch$,
\[
  \RGamma : \DX \longrightarrow \D(\Mod\,k),
\]
(which equals the derived direct image functor $\Rf$) has a
right-adjoint 
\[
  f^! : \D(\Mod\,k) \longrightarrow \DX.
\]
It is easy to see $f^!k \cong \omega[d]$; cf.\ \cite[rmk.\
5.5]{NeemanDuality}.  Observe also that since the singularities of
$\sch$ are Gorenstein, $\omega$ is an invertible sheaf.

Now consider the following sequence of natural isomorphisms which is
taken from \cite[sec.\ 4]{KrauAR2}.  The object $\testcomp$ is in
$\DX$. 
\begin{eqnarray*}
    & & \Hom_{\DsmallX}(\rightobj,\testcomp)^{\prime} \\
    & & \;\;\;\;\; \cong
        \Hom_{\Dsmall(\Modsmall\,k)}(
          \RGamma(\RSHom_{\sch}(\rightobj,\testcomp)),k) \\
    & & \;\;\;\;\; \cong
        \Hom_{\DsmallX}(
          \RSHom_{\sch}(\rightobj,\testcomp),f^!k) \\
    & & \;\;\;\;\; \stackrel{\rm (a)}{\cong}
        \Hom_{\DsmallX}(
          \RSHom_{\sch}(
            \rightobj,\cO_{\sch}) \LTensor_{\sch} \testcomp,f^!k) \\
    & & \;\;\;\;\; \stackrel{\rm (b)}{\cong}
        \Hom_{\DsmallX}(
          \testcomp,\RSHom_{\sch}(
            \RSHom_{\sch}(\rightobj,\cO_{\sch}),f^!k)) \\
    & & \;\;\;\;\; \stackrel{\rm (c)}{\cong}
        \Hom_{\DsmallX}(\testcomp,\rightobj \LTensor_{\sch} f^!k). \\
\end{eqnarray*}
Here (b) is by adjointness between $\LTensor$ and $\RSHom$, while (a)
and (c) hold because they clearly hold if $\rightobj$ is a locally
free coherent sheaf, and therefore also hold for the given $\rightobj$
because it has a bounded resolution of locally free coherent sheaves,
so is finitely built in $\DX$ from such sheaves.

Now, $f^!k$ is $\omega[d]$ and $\omega$ is an invertible sheaf.
Hence $\rightobj \LTensor_{\sch} f^!k$ is just $(\rightobj
\otimessch \omega)[d]$.  Inserting this along with
$\testcomp = \testobj[i]$ with $\testobj$ in $\QCoh(\sch)$ gives
natural isomorphisms
\[
  \Ext_{\QCohsmall\,\sch}^i(\rightobj,\testobj)^{\prime}
  \cong \Ext_{\QCohsmall\,\sch}^{d-i}(
          \testobj,\rightobj \otimessch \omega),
\]
proving the lemma.
\end{proof}

\begin{Theorem}
\label{thm:ARs}
Let $\sch$ be a smooth projective scheme of dimension $d \geq 1$ over
the field $k$.  Let $\omega$ be the dualizing sheaf.

Let $\rightobj$ be an indecomposable coherent sheaf.  Then there is an
\ARs\ in the category of quasi-coherent sheaves,
\[
  0 \rightarrow (\Sigma^{d-1}\rightobj) \otimessch \omega
    \longrightarrow \midobj
    \longrightarrow \rightobj
    \rightarrow 0.
\]
Here $\Sigma^{d-1}\rightobj$ is the $(d-1)$'st syzygy in a minimal
injective resolution of $\rightobj$ in the category of quasi-coherent
sheaves.
\end{Theorem}

\begin{proof}
The proof will give slightly more than stated: Let $\sch$ be a
projective scheme over $k$ of dimension $d \geq 1$ with Gorenstein
singularities, and let $\rightobj$ be an indecomposable coherent sheaf
which has a bounded resolution of locally free coherent sheaves.  Then
I will prove that the indicated \ARs\ exists.  (In the smooth case,
each coherent sheaf has a resolution as required by
\cite[exer.\ III.6.5]{Hartshorne}).

The category of quasi-coherent sheaves $\QCoh(\sch)$ is clearly
$k$-linear, and by \cite[lem.\ 3.1]{Illusie}, it is a
Gro\-then\-di\-eck category, that is, a cocomplete abelian category
with a generator and exact filtered colimits.  Hence it has injective
envelopes by \cite[thm.\ 10.10]{Popescu}, and in particular, it has
enough injectives.

Also, the category of coherent sheaves is an abelian category so has
split idempotents.  It also has finite dimensional $\Hom$ sets, as
follows e.g.\ from Serre finiteness, \cite[thm.\
III.5.2(a)]{Hartshorne}.  Hence the endomorphism ring of the
indecomposable sheaf $\rightobj$ is local; cf.\ \cite[p.\
52]{RingelTame}. 

Finally, if I set $i = 0$ in lemma \ref{lem:duality}, then I get the
natural equivalence
\[
  \HomQCoh(\rightobj,-)^{\prime}
  \simeq 
  \ExtQCoh^d(-,\rightobj \otimessch \omega).
\]
This is equation \eqref{equ:equivalence} of proposition \ref{pro:ARs},
with $\leftobj = \rightobj \otimessch \omega$.  

So all the conditions of proposition \ref{pro:ARs} are satisfied, and
hence, there is a short exact sequence in $\QCoh(\sch)$,
\[
  0 \rightarrow \Sigma^{d-1}(\rightobj \otimessch \omega)
    \stackrel{\leftobjmor}{\longrightarrow} \midobj
    \stackrel{\midobjmor}{\longrightarrow}  \rightobj
    \rightarrow 0,
\]
where $\midobjmor$ is right almost split.  If I construct
$\Sigma^{d-1}(\rightobj \otimessch \omega)$ by means of a minimal
injective resolution of $\rightobj \otimessch \omega$, then this
sequence equals the one in the theorem because $\omega$ is invertible.
So to finish the proof, I must show that $\leftobjmor$ is left almost
split.  Since $\midobjmor$ is right almost split, it is enough to show
that $\EndQCoh(\Sigma^{d-1}(\rightobj \otimessch \omega))$ is a local
ring, by classical \ARth.

For this, note that the minimal injective resolution gives short exact
sequences
\[
  0 \rightarrow \Sigma^{\ell}(\rightobj \otimessch \omega)
    \longrightarrow \injobj^{\ell}
    \longrightarrow \Sigma^{\ell+1}(\rightobj \otimessch \omega)
    \rightarrow 0
\]
for $\ell = 0, \ldots, d-2$, where $\Sigma^{\ell}(\rightobj
\otimessch \omega) \longrightarrow \injobj^{\ell}$ is an injective
envelope.  Hence, starting with the knowledge that
$\EndQCoh(\rightobj \otimessch \omega) \cong
\EndQCoh(\rightobj)$ is local, successive uses of lemma
\ref{lem:End} will prove that all the rings
\[
  \EndQCoh(\rightobj \otimessch \omega), 
  \ldots,
  \EndQCoh(\Sigma^{d-1}(\rightobj \otimessch \omega)) 
\]
are local, provided that I can show
\[
  \ExtQCoh^1(
    \injobj^{\ell},\Sigma^{\ell}(\rightobj \otimessch \omega)) = 0
\]
for $\ell = 0, \ldots, d-2$, that is,
\begin{equation}
\label{equ:Ext_vanishing}
  \ExtQCoh^{\ell+1}(\injobj^{\ell},\rightobj \otimessch \omega) = 0
  \; \mbox{ for } \; \ell = 0, \ldots, d-2.
\end{equation}
However, if $\injobj$ is any injective, then lemma \ref{lem:duality}
gives
\[
  \ExtQCoh^j(\injobj,\rightobj \otimessch \omega) 
  \cong \ExtQCoh^{d-j}(\rightobj,\injobj)^{\prime} = 0
\]
for $j = 1, \ldots, d-1$, and this implies \eqref{equ:Ext_vanishing}.
\end{proof}

Theorem \ref{thm:ARs} allows me to recover the following result
known from \cite{ReitenVandenBergh}.

\begin{Corollary}
\label{cor:ARs}
Let $\sch$ be a smooth projective curve over the field $k$.  Then
the category of coherent sheaves has right and left \ARs s.
\end{Corollary}

\noindent
Indeed, if $\rightobj$ and $\leftobj$ are indecomposable
coherent sheaves, then there are \ARs s of coherent sheaves,
\[
  0 \rightarrow \rightobj \otimessch \omega
    \longrightarrow \midobj
    \longrightarrow \rightobj
    \rightarrow 0
\]
and
\[
  0 \rightarrow \leftobj 
    \longrightarrow \midobjprime
    \longrightarrow \leftobj \otimessch \omega^{-1}
    \rightarrow 0.
\]

\noindent
{\bf Acknowledgement.}  I would like to thank Michel Van den Bergh for
several comments to a preliminary version, and Apostolos Beligiannis
and Henning Krause for communicating \cite{Beligiannis} and
\cite{KrauAR2}.

The diagrams were typeset with Paul Taylor's {\tt diagrams.tex}.

\end{document}